\documentclass[a4paper,12pt]{amsart}

\pdfoutput = 1 % sets the output file to be in PDF format (instead of DVI)

\usepackage[
a4paper ,
top    = 3.3cm ,
bottom = 2cm ,
left   = 2.5cm ,
right  = 2.5cm ,
]{geometry}

\usepackage[T1]{fontenc} % allows to type accented characters directly into the .tex file

\usepackage[american]{babel}
\usepackage{mathtools} % amsmath improvement
\mathtoolsset{mathic} % math italics correction

\usepackage[ScaleRM = 1.06]{libertinus} % serif text font
\usepackage[scale = 0.99]{FiraSans} % sans serif text font
\usepackage[scale = 1]{inconsolata} % sans serif monospaced ("typewriter") font

\usepackage{eulerpx} % math font
\usepackage{bm} % improved bold math

\makeatletter
\g@addto@macro\bfseries{\boldmath}
\makeatother

\allowdisplaybreaks[3] % allow page breaks in multiline equations

\makeatletter
\renewcommand\operator@font{\sf}
\makeatother

\DeclareTextFontCommand{\noun}{\scshape} % small caps for proper nouns

\newcommand{\mybar}[3]{
\mathrlap{\hspace{#2}\overline{\scalebox{#1}[1]{\phantom{\ensuremath{#3}}}}}\ensuremath{#3}
}
\newcommand{\Abar}{ \mybar{0.7}{1.5pt}{A} }

\newcommand{\EE}{ \varmathbb{E} }
\newcommand{\KK}{ \varmathbb{K} }
\newcommand{\CC}{ \varmathbb{C} }
\newcommand{\RR}{ \varmathbb{R} }

\newcommand{\LL}    { \mathsf{L} }
\newcommand{\St}    { \mathsf{S} }
\newcommand{\Cotype}{ \mathsf{C} }
\newcommand{\Kah}   { \mathsf{K} }

\bmdefine{\aa}{ \mathrm{a} }
\bmdefine{\ee}{ \mathrm{e} }
\bmdefine{\tt}{ \mathrm{t} }
\bmdefine{\xx}{ \mathrm{x} }
\bmdefine{\ww}{ \mathrm{w} }
\bmdefine{\yy}{ \mathrm{y} }
\bmdefine{\zz}{ \mathrm{z} }

\newcommand{\Aab}{ \norm{A}_{ a , b } }

\DeclarePairedDelimiter{\abs}{\lvert}{\rvert}
\DeclarePairedDelimiter{\norm}{\lVert}{\rVert}

\input{hksqrt} % root symbol with hook

\numberwithin{equation}{section} % numbering of equations includes section number

\usepackage[pagebackref]{hyperref}
\hypersetup{
allcolors     = blue ,
breaklinks    = true ,
colorlinks    = true ,
pdfencoding   = unicode ,
pdfpagelayout = OneColumn ,
pdfstartview  = {FitH} ,
}

\usepackage{bookmark} % PDF bookmarks

\usepackage[
numeric ,
msc-links ,
nobysame ,
]{amsrefs}

\makeatletter
\def\print@backrefs#1{\space\SentenceSpace[cited on page(s) \csname br@#1\endcsname]}
\makeatother

\BibSpec{book}{
    +{}  {\PrintPrimary}                {transition}
    +{,} { \textit}                     {title}
    +{.} { }                            {part}
    +{:} { \textit}                     {subtitle}
    +{,} { \PrintEdition}               {edition}
    +{}  { \PrintEditorsB}              {editor}
    +{,} { \PrintTranslatorsC}          {translator}
    +{,} { \PrintContributions}         {contribution}
    +{,} { }                            {series}
    +{,} { \voltext}                    {volume}
    +{,} { }                            {publisher}
    +{,} { }                            {organization}
    +{,} { }                            {address}
    +{,} { \PrintDateB}                 {date}
    +{,} { }                            {status}
    +{,} { \PrintDOI}                   {doi}
    +{}  { \parenthesize}               {language}
    +{}  { \PrintTranslation}           {translation}
    +{;} { \PrintReprint}               {reprint}
    +{.} { }                            {note}
    +{.} {}                             {transition}
    +{}  {\SentenceSpace \PrintReviews} {review}
}

\usepackage[doipre={}]{uri} % DOI hyperlinks for bibliography items

\usepackage{enumitem} % customized lists

\usepackage[framemethod = tikz]{mdframed}
\mdfdefinestyle{barra}{
backgroundcolor   = white ,
linewidth         = 3pt ,
topline           = false ,
bottomline        = false ,
rightline         = false ,
innertopmargin    = -6pt ,
innerbottommargin = 0pt ,
innerrightmargin  = 0pt ,
innerleftmargin   = 5pt ,
}

\theoremstyle{plain}
\newmdtheoremenv[ default , style = barra ]{theorem}[dummy]{Theorem}
\newmdtheoremenv[ default , style = barra ]{corollary}[dummy]{Corollary}
\newmdtheoremenv[ default , style = barra ]{proposition}[dummy]{Proposition}
\newmdtheoremenv[ default , style = barra ]{lemma}[dummy]{Lemma}

\theoremstyle{definition}
\newmdtheoremenv[ default , style = barra ]{definition}[dummy]{Definition}
\newmdtheoremenv[ default , style = barra ]{example}[dummy]{Example}
\newmdtheoremenv[ default , style = barra ]{remark}[dummy]{Remark}
\newmdtheoremenv[ default , style = barra ]{problem}{Problem}

\usepackage{calc} % allows to use infix notation arithmetic in arguments to macros
\newlength{\figwidth}
\newlength{\figwidthwithoutcaption}
\begin{document}

\title[About the real anisotropic \noun{Littlewood}'s $4 / 3$ inequality]{The sharp constants in the real anisotropic \\ \noun{Littlewood}'s $4 / 3$ inequality and applications}

\author[N.~Caro-Montoya]{Nicolás Caro-Montoya}
\address{Departamento de Matemática \\ Centro de Ciências Exatas e da Natureza \\ Universidade Federal de Pernambuco \\ Avenida Jornalista Aníbal Fernandes, S/N - Cidade Universitária \\ Recife/PE - Brasil - 50740-560}
\email[N.~Caro-Montoya]{jorge.caro@ufpe.br}

\author[D.~Núñez-Alarcón]{Daniel Núñez-Alarcón}
\address{Departamento de Matemáticas \\ Universidad Nacional de Colombia \\ Bogotá, Colombia - 111321}
\email[D.~Núñez-Alarcón]{dnuneza@unal.edu.co}

\author[D.~Serrano-Rodríguez]{Diana Serrano-Rodríguez}
\email[D.~Serrano-Rodríguez]{diserranor@unal.edu.co}

\dedicatory{Dedicated to the memory of \,\noun{Jhazaira Mantilla~Pérez}}


\begin{abstract}

The real anisotropic \noun{Littlewood}'s $4 / 3$ inequality is an extension of a famous result obtained in 1930 by \noun{J.~E.~Littlewood}. It asserts that, for $a , b \in ( 0 , \infty )$, the following conditions are equivalent:

\medskip

\begin{itemize}[itemsep = 8pt]

\item There is an optimal constant $\LL_{ a , b }^{ \RR } \in [ 1 , \infty )$ such that
\[
\Biggl ( \, \sum_{ k = 1 }^{ \infty } \biggl ( \, \sum_{ j = 1 }^{ \infty } \, \abs[\big]{ A \bigl ( \ee^{ (k) } , \ee^{ (j) } \bigr ) }^a \biggr )^{ \frac{b}{a} } \Biggr )^{ \frac{1}{b} } \leq \LL_{ a , b }^{ \RR } \cdot \norm{A}
\]
for every continuous bilinear form $A \colon c_0 \times c_0 \to \RR$.

\item The values $a , b$ satisfy $a , b \geq 1$ and $\frac{1}{a} + \frac{1}{b} \leq \frac{3}{2}$\,.

\end{itemize}
\bigskip
Several authors have obtained the values of $\LL_{ a , b }^{ \RR }$ for diverse pairs $( a , b )$. In this paper we provide the complete list of such optimal values, as well as new estimates for $\LL_{ a , b }^{ \CC }$ (the analog for continuous $\CC$-bilinear forms), which are exact in several cases. As an application we prove, in terms of the values $\LL_{ 1 , r }^{ \CC }$\,, a variant of \noun{Khinchin}'s inequality for \noun{Steinhaus} variables, and we provide estimates for the $( q , s )$-cotype constants of the spaces $\ell_1 ( \KK )$ (with $\KK = \RR$ or $\CC$) in terms of the values $\LL_{ 1 , q }^{ \RR }$\,.

\end{abstract}

\maketitle % displays title and abstract

\tableofcontents

\section{Introduction}

\noindent The origins of the theory of absolutely summing multilinear operators are linked to the so-called \noun{Littlewood}'s $4 / 3$ inequality~\cite{l1930}; for a detailed introduction to the theory of absolutely summing operators see~\cite{diestel}. On the other hand, the multilinear theory of absolutely summing operators has been recently explored in different contexts by several authors (see~\cites{bayart2, advances, bps, pilar, bot, psst} and the references therein), with unexpected applications in other fields such as quantum information theory (more precisely, quantum XOR games~\cite{monta}), computational learning theory~\cites{ACKSW, ADEP, KSVZ} and results related to the study of \noun{Dirichlet} series (more precisely, questions concerning the $n$-dimensional \noun{Bohr} radius; see~\cite{bps}). The study of sharp constants for inequalities of the type \noun{Littlewood}'s $4 / 3$ (and its multilinear version, namely the \noun{Bohnenblust--Hille}'s inequality), among other techniques, allowed to achieve such applications. Driven by this, during the last decade several works have emerged dealing with generalizations of \noun{Littlewood}'s $4 / 3$ inequality to the multilinear context and also to other sequence spaces (see~\cites{nusete,rapser} and the references therein).

All along this paper $\KK$ stands for either $\RR$ or $\CC$, and $c_0 = c_0 ( \KK )$, the \noun{Banach} space of null sequences. By convention, when $s = \infty$, the value $\bigl ( \sum_j \, \abs{ x_j }^s \bigr )^{ 1 / s }$ represents the supremum of the values $\abs{ x_j }$. We also define $\frac{1}{0} = \infty$ and $\frac{1}{ \infty } = 0$, and for $s \in [ 1 , \infty ]$ we denote by $s^{ \ast }$ the conjugate index of $s$, so that $\frac{1}{s} + \frac{1}{ s^{ \ast } } = 1$. Finally, we denote the sequence of canonical vectors of $c_0$ by $\bigl ( \ee^{ (n) } \bigr )_{ n \geq 1 }$\,.

\noun{Littlewood}'s $4 / 3$ inequality~\cite{l1930} asserts that there is a constant $\LL_{ \KK } \in [ 1 , \infty )$ such that
\[
\Biggl ( \, \sum_{ k , j = 1 }^{ \infty } \, \abs[\big]{ A \bigl ( \ee^{ (k) } , \ee^{ (j) } \bigr ) }^{ \frac{4}{3} } \Biggr )^{ \frac{3}{4} } \leq \LL_{ \KK } \cdot \norm{A}
\]
for every continuous bilinear form $A \colon c_0 \times c_0 \to \KK$, where $\norm{ \cdot }$ is the usual norm given by
\[
\norm{A} = \sup \, \bigl \{ \abs{ A( \xx \, , \yy ) } \colon \xx , \yy \in c_0 \ \ ; \ \ \norm{ \xx }_{ c_0 } \, , \norm{ \yy }_{ c_0 } \leq 1 \bigr \} \, .
\]
It is well known that the exponent $4 / 3$ is optimal, and it is shown in~\cite{diniz} that the value $\LL_{ \RR } = \sqrt{2}$ is also optimal. For complex scalars we only know that the corresponding optimal value satisfies $\LL_{ \CC } \leq \frac{2}{ \sqrt{ \pi } }$\,.

However, the optimality feature of the exponent $4 / 3$ is related to the specific configuration of this inequality. A key issue in the theory has been to investigate optimality ranges for summing inequalities involving the so-called \emph{anisotropic exponents}. More specifically, our objective is to control the quantity
\begin{equation} \label{normaAab}
\Biggl ( \, \sum_{ k = 1 }^{ \infty } \biggl ( \, \sum_{ j = 1 }^{ \infty } \, \abs[\big]{ A \bigl ( \ee^{ (k) } , \ee^{ (j) } \bigr ) }^a \biggr )^{ \frac{b}{a} } \Biggr )^{ \frac{1}{b} } \coloneq \Biggl ( \, \sum_{ k = 1 }^{ \infty } \Biggl [ \biggl ( \, \sum_{ j = 1 }^{ \infty } \, \abs[\big]{ A \bigl ( \ee^{ (k) } , \ee^{ (j) } \bigr ) }^a \biggr )^{ \frac{1}{a} } \Biggr ]^b \, \Biggr )^{ \frac{1}{b} } \tag{$\maltese$}
\end{equation}
for all norm-$1$ continuous bilinear form $A \colon c_0 \times c_0 \to \KK$. This problem has been addressed through several approaches permeating the theory, with significant advances obtained over the past decade.

The following anisotropic \noun{Littlewood}'s $4 / 3$ inequality is a particular case of~\cite{psst}*{Theorem~5.1}:

\begin{proposition} \label{Lab-weak}

Given \( a , b \in ( 0 , \infty ) \), the following assertions are equivalent\textup{:}

\begin{enumerate}[
itemsep     = 8pt ,
labelindent = 3mm ,
leftmargin  = * ,
label       = \textup{\texttt{\alph*)}} ,
align       = left ,
partopsep   = 2mm ,
]

\item There is a constant \( \LL \in [ 1 , \infty ) \) such that
\[
\Biggl ( \, \sum_{ k = 1 }^{ \infty } \biggl ( \, \sum_{ j = 1 }^{ \infty } \, \abs[\big]{ A \bigl ( \ee^{ (k) } , \ee^{ (j) } \bigr ) }^a \biggr )^{ \frac{b}{a} } \Biggr )^{ \frac{1}{b} } \leq \LL \cdot \norm{A}
\]
for every continuous bilinear form \( A \colon c_0 \times c_0 \to \KK \).

\item The values \( a \) and \( b \) satisfy \( a , b \geq 1 \) and \( \frac{1}{a} + \frac{1}{b} \leq \frac{3}{2} \)\,.

\end{enumerate}

\end{proposition}

\noindent This result improves that of~\cite{ABPS}*{Theorem~1.1}, where the same equivalence is proved under the additional assumption $a , b \in [ 1 , 2 ]$. Observe that, by taking $a = b = 4 / 3$ in Proposition~\ref{Lab-weak}, one recovers \noun{Littlewood}’s $4 / 3$ inequality.

From now on, for $a , b \in ( 0 , \infty ]$ we denote the value~\eqref{normaAab} by $\Aab \in [ 0 , \infty ]$. In the particular case $a , b \geq 1$, the set of $A$ with $\Aab < \infty$ is a vector space with norm $\norm{ \cdot }_{ a , b }$ (see~\cite{BePa61}*{p.~302}). Moreover, by our convention, for $a , b \in ( 0 , \infty )$ we have
\begin{align*}
\norm{A}_{ a , \infty } & = \sup_{ k \geq 1 } \, \Biggl ( \, \sum_{ j = 1 }^{ \infty } \abs[\big]{ A \bigl ( \ee^{ (k) } , \ee^{ (j) } \bigr ) }^a \Biggr )^{ \frac{1}{a} } \, ; \\[2mm]
\norm{A}_{ \infty , b } & = \Biggl ( \, \sum_{ k = 1 }^{ \infty } \biggl ( \, \sup_{ j \geq 1 } \ \abs[\big]{ A \bigl ( \ee^{ (k) } , \ee^{ (j) } \bigr ) } \biggr )^b \, \Biggr )^{ \frac{1}{b} } \, ,
\intertext{and}
\norm{A}_{ \infty , \infty } & = \sup_{ k , j \geq 1 } \ \abs[\big]{ A \bigl ( \ee^{ (k) } , \ee^{ (j) } \bigr ) } \, .
\end{align*}

\begin{remark} \label{0<s<1;(s,infty)and(infty,s)nonadmissible}

By~\cite{zalduendo}*{Proposition~1} there exists $\varphi \in c_0'$ such that $\sum_{ n = 1 }^{ \infty } \, \abs[\big]{ \varphi \bigl ( \ee^{ (n) } \bigr ) }^s$ diverges for all $s \in ( 0 , 1 )$. If $A_1 , A_2 \colon c_0 \times c_0 \to \KK$ are defined by
\begin{align*}
A_1 ( \xx \, , \yy ) & = \varphi ( \yy ) x_1 \, ; \\
A_2 ( \xx \, , \yy ) & = \varphi ( \xx ) y_1 \, ,
\end{align*}
then it is straightforward to check that $A_1$ and $A_2$ are continuous bilinear forms with $\norm{ A_1 } = \norm{ A_2 } = \norm{ \varphi }_{ c_0' }$\,; moreover, for all $s \in ( 0 , 1 )$ we have
\[
\norm{ A_1 }_{ s , \infty } = \norm{ A_2 }_{ \infty , s } = \sum_{ n = 1 }^{ \infty } \, \abs[\big]{ \varphi \bigl ( \ee^{ (n) } \bigr ) }^s = \infty \, .
\]
Consequently, for all $s \in ( 0 , 1 )$ the pairs $( a , b ) = ( s , \infty )$ and $( a , b ) = ( \infty , s )$ do not satisfy condition \texttt{a)} in Proposition~\ref{Lab-weak}.

\end{remark}

\noindent Given $a , b \in ( 0 , \infty ]$ satisfying condition \texttt{a)} in Proposition~\ref{Lab-weak}, we denote the corresponding optimal value of $\LL$ by $\LL_{ a , b }^{ \KK } \in [ 1 , \infty )$. In the real case we have~\cite{ABPS}*{Theorem~6.3}:
\begin{equation} \label{Lab^R-incomplete}
\LL_{ a , b }^{ \RR } = 2^{ \frac{1}{a} + \frac{1}{b} - 1 } , \quad \textnormal{for} \ a , b \in [ 1 , 2 ] \ \ \textnormal{with} \ \ \frac{1}{a} + \frac{1}{b} \leq \frac{3}{2} \, .
\end{equation}
On the other hand, for complex scalars, the determination of the exact values of the optimal constants involved is probably a difficult task. In this case the only known optimal values are~\cite{blei}*{p.~31}
\begin{equation} \label{L_12^C;L_21^C}
\LL_{ 1 , 2 }^{ \CC } = \LL_{ 2 , 1 }^{ \CC } = \frac{2}{ \sqrt{ \pi } } \, ,
\end{equation}
and in general, for $a , b \in [ 1 , 2 ]$ with $\frac{1}{a} + \frac{1}{b} \leq \frac{3}{2}$\,, we have the estimate~\cite{ABPS}*{Remark 6.4}
\begin{equation} \label{Lab^C}
\LL_{ a , b }^{ \CC } \leq \biggl ( \frac{4}{ \pi } \biggr )^{ \frac{1}{a} + \frac{1}{b} - 1 } \, . \pagebreak
\end{equation}

\begin{figure}[tbh!]
%
% sets the figure centered with respect to the graphic without the caption, which is precisely the width of Figure 3 (which has no caption)
\settowidth{\figwidth}{\includegraphics{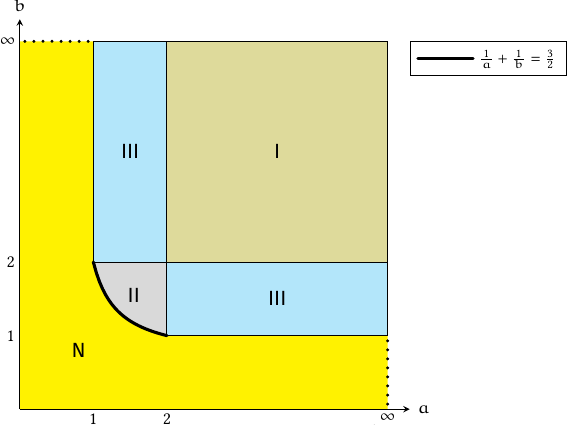}}
\hspace*{\figwidth - \figwidthwithoutcaption}
\includegraphics{Figure1.pdf}
\caption{Region~\textsf{N} (including the dotted open segments but not other points of its boundary) corresponds to non-admissible exponents in the real and complex anisotropic \noun{Littlewood}'s $4 / 3$ inequality (Proposition~\ref{Lab-weak} and Remark~\ref{0<s<1;(s,infty)and(infty,s)nonadmissible}). Region~\textsf{I} (closed) corresponds to Remark~\ref{Lab^K=1;ab>=2}; in this region the optimal constant, both in the real and complex case, has always value $1$. Region~\textsf{II} (closed) corresponds to Equation~\eqref{Lab^R-incomplete}, so that the value of the optimal constant in the real case is $2^{ \frac{1}{a} + \frac{1}{b} - 1 }$. Finally, Region~\textsf{III} corresponds to Problem~\ref{Problem1}, whereas the union of Regions \textsf{II} and \textsf{III} minus the two points $( 1 , 2 )$ and $( 2 , 1 )$ corresponds to Problem~\ref{Problem2}; see Equation~\eqref{L_12^C;L_21^C}.}
\label{fig1}
\end{figure}

\begin{remark} \label{Lab^K=1;ab>=2}

From Equations~\eqref{Lab^R-incomplete} and~\eqref{Lab^C} we get $\LL_{ 2 , 2 }^{ \KK } = 1$. Therefore, if $a , b \in [ 2 , \infty ]$, then from the embedding of the $\ell_p$ spaces we conclude that $\LL_{ a , b }^{ \KK } = 1$ (see Figure~\ref{fig1}).

\end{remark}

\begin{remark} \label{minkowski}

By \noun{Minkowski}'s inequality~\cite{garling}*{Corollary~5.4.2}, for any $a , b \in ( 0 , \infty )$ with $a \leq b$ and any double sequence $( x_{ k , j } )_{ k , j \geq 1 }$ of elements in $\KK$ we have
\[
\Biggl ( \, \sum_{ k = 1 }^{ \infty } \biggl ( \, \sum_{ j = 1 }^{ \infty } \, \abs{ x_{ k , j } }^a \biggr )^{ \frac{b}{a} } \Biggr )^{ \frac{1}{b} } \leq \Biggl ( \, \sum_{ k = 1 }^{ \infty } \biggl ( \, \sum_{ j = 1 }^{ \infty } \, \abs{ x_{ j , k } }^b \biggr )^{ \frac{a}{b} } \Biggr )^{ \frac{1}{a} } \, .
\]
Consequently, given a continuous bilinear form $A \colon c_0 \times c_0 \to \KK$, the corresponding transpose map $\Abar$ is also continuous (with $\norm{ \Abar } = \norm{A}$) and satisfies $\Aab \leq \norm{ \Abar }_{ b , a }$\,. In this way, if $\LL_{ b , a }^{ \KK }$ exists, then $\LL_{ a , b }^{ \KK }$ also exists and satisfies $\LL_{ a , b }^{ \KK } \leq \LL_{ b , a }^{ \KK }$\,.

\end{remark}
\bigskip
\noindent The facts above lead us to pose the following questions:

\begin{problem}[see Figure~\ref{fig1}] \label{Problem1}

What are the values of $\LL_{ a , b }^{ \RR }$ for $a , b \in [ 1 , \infty ]$ with $a < 2 < b$ or $b < 2 < a$?

\end{problem}

\begin{problem} \label{Problem2}

What are the values of $\LL_{ a , b }^{ \CC }$ for $a , b \in [ 1 , \infty ]$ with $\frac{1}{a} + \frac{1}{b} \leq \frac{3}{2}$?

\end{problem}

\pagebreak

\noindent In this paper we solve Problem~\ref{Problem1}, thus providing a complete description of the best constants in the real anisotropic \noun{Littlewood}'s $4 / 3$ inequality:

\begin{theorem} \label{main-Lab^R}

Let \( a , b \in [ 1 , \infty ] \) with \( \frac{1}{a} + \frac{1}{b} \leq \frac{3}{2} \)\,. For every continuous bilinear form \( A \colon c_0 \times c_0 \to \RR \) we have
\[
\Aab \leq 2^{ \max \{ 0 , \frac{1}{a} + \frac{1}{b} - 1 \} } \cdot \norm{A} \, .
\]
Moreover, the value \( 2^{ \max \{ 0 , \frac{1}{a} + \frac{1}{b} - 1 \} } \) is optimal.

\end{theorem}

\noindent This paper is organized as follows: Section~2 is devoted to the proof of our main result (Theorem~\ref{main-Lab^R}), which uses the variant of \noun{Khinchin}'s inequality from~\cite{pelsansan} together with an interpolation result (Proposition~\ref{interpol}): this provides a complete description of the best constants in the real anisotropic \noun{Littlewood}'s $4 / 3$ inequality, thus completely solving Problem~\ref{Problem1}; in addition, some of these ideas are adapted to the setting of complex scalars, which allows us to partially solve Problem~\ref{Problem2}. In Section~3 we prove a variant (in the sense of~\cite{pelsansan}) of \noun{Khinchin}'s inequality for \noun{Steinhaus} variables, involving the best constants previously discussed. Finally, we show in Section~4 how the optimal $( q , s )$-cotype constants of the space $\ell_{1} ( \KK )$ can be estimated in terms of the best constants from Section~2.

\section{The main result}

\noindent The first part of this section is devoted to the proof of Theorem~\ref{main-Lab^R}. The auxiliary results required are stated and proved for arbitrary (real or complex) scalars; among these, Lemma~\ref{La1^K<=2^1/a} involves a constant that can be improved in the specialized case $\KK = \CC$. This in turn allows us to obtain, for complex scalars, a partial counterpart of the main result (Section~\ref{subsection-complex}).

We need the following variant of \noun{Khinchin}'s inequality, which is a particular case of~\cite{pelsansan}*{Lemma~1 and Proposition~1}. Its proof strongly uses the so-called Contraction Principle from~\cite{diestel}*{p.~231}.

\begin{proposition} \label{radem-cota-2^1/a}

Given \( a \in [ 2 , \infty ] \) we have, for any complex sequence \( ( c_n )_{ n \geq 1 } \) and all \( N \geq 1 \):
\begin{equation} \label{eq-radem-cota-2^1/a}
\Biggl ( \, \sum_{ n = 1 }^N \, \abs{ c_n }^a \Biggr )^{ \frac{1}{a} } \leq 2^{ \frac{1}{a} } \int_0^1 \abs[\Bigg]{ \sum_{ n = 1 }^N r_n (t) \,c_n } \, \mathrm{d} t \, .
\end{equation}
Moreover, the value \( 2^{ \frac{1}{a} } \) is optimal.

\end{proposition}

\noindent Above, as usual, $( r_n \colon [ 0 , 1 ] \to \RR )_{ n \geq 1 }$ is the sequence of independent and identically distributed random variables defined by
\begin{equation} \label{rademacher}
r_n (t) = \operatorname{sign} \bigl ( \sin ( 2^n \pi t ) \bigr ) \, ,
\end{equation}
the so-called \emph{\noun{Rademacher} functions}.

As a first step in the proof of our main result, we deal with the following boundary case:

\begin{lemma} \label{La1^K<=2^1/a}

For all \( a \in [ 2 , \infty ] \) and every continuous bilinear form \( A \colon c_0 \times c_0 \to \KK \) we have
\[
\norm{A}_{ a , 1 } = \sum_{ k = 1 }^{ \infty } \Biggl ( \, \sum_{ j = 1 }^{ \infty } \, \abs[\big]{ A \bigl ( \ee^{ (k) } , \ee^{ (j) } \bigr ) }^a \Biggr )^{ \frac{1}{a} } \leq 2^{ \frac{1}{a} } \cdot \norm{A} \, .
\]

\end{lemma}

\begin{proof}

Let $N$ be a positive integer. From inequality~\eqref{eq-radem-cota-2^1/a} we get, for all $a \in [ 2 , \infty ]$:
\[
\sum_{ k = 1 }^N \Biggl ( \, \sum_{ j = 1 }^N \, \abs[\big]{ A \bigl ( \ee^{ (k) } , \ee^{ (j) } \bigr ) }^a \Biggr )^{ \frac{1}{a} } \leq 2^{ \frac{1}{a} } \sum_{ k = 1 }^N \int_0^1 \abs[\Bigg]{ \sum_{ j = 1 }^N r_j (t) \, A \bigl ( \ee^{ (k) } , \ee^{ (j) } \bigr ) } \, \mathrm{d} t \, .
\]
We claim that the right side term above is bounded above by $2^{ \frac{1}{a} } \cdot \norm{A}$; having this, it is enough to take $N \to \infty$ on the left side to finish the proof.

To prove our claim, recall that for every $\varphi \in c_0'$ we have $\sum_{ k = 1 }^{ \infty } \abs[\big]{ \varphi \bigl ( \ee^{ (k) } \bigr ) } \leq \norm{ \varphi }_{ c_0' }$\,. Therefore we have

\begin{align} \label{bound-integrand}
2^{ \frac{1}{a} } \sum_{ k = 1 }^N \int_0^1 \abs[\Bigg]{ \sum_{ j = 1 }^N r_j (t) \, A \bigl ( \ee^{ (k) } , \ee^{ (j) } \bigr ) } \, \mathrm{d} t & = 2^{ \frac{1}{a} } \int_0^1 \sum_{ k = 1 }^N \, \abs[\Bigg]{ A \Biggl ( \ee^{ (k) } , \sum_{ j = 1 }^N r_j (t) \cdot \ee^{ (j) } \Biggr ) } \, \mathrm{d} t \notag \\[2mm]
& \leq 2^{ \frac{1}{a} } \int_0^1 \norm[\Bigg]{ A \Biggl ( \ \bm{ \relbar } \ , \sum_{ j = 1 }^N r_j (t) \cdot \ee^{ (j) } \Biggr ) }_{ c_0' } \! \mathrm{d} t \, . \tag{$\ddagger$}
\end{align}
Finally, for any $t \in [ 0 , 1 ]$ we have $\norm[\big]{ \sum_{ j = 1 }^N r_j (t) \cdot \ee^{ (j) } }_{ c_0 } \leq 1$, and thus the integrand in the upper bound~\eqref{bound-integrand} satisfies
\begin{align*}
\norm[\Bigg]{ A \Biggl ( \ \bm{ \relbar } \ , \sum_{ j = 1 }^N r_j (t) \cdot \ee^{ (j) } \Biggr ) }_{ c_0' } & = \sup_{ \norm{ \xx }_{ c_0 } \leq 1 } \ \abs[\Bigg]{ A \Biggl ( \xx \, , \sum_{ j = 1 }^N r_j (t) \cdot \ee^{ (j) } \Biggr ) } \\[2mm]
& \leq \norm{A} \, . \qedhere
\end{align*}

\end{proof}

\begin{remark} \label{L1infty=Linfty1=1}

Taking $a = \infty$ in Lemma~\ref{La1^K<=2^1/a} yields $\LL_{ \infty , 1 }^{ \KK } = 1$, whereas the chain of inequalities
\[
\sup_{ k \geq 1 } \, \sum_{ j = 1 }^{ \infty } \,\abs[\big]{ A \bigl ( \ee^{ (k) } , \ee^{ (j) } \bigr ) } \leq \sup_{ k \geq 1 } \ \norm[\big]{ A \bigl ( \ee^{ (k) } , \bm{ \relbar } \, \bigr ) }_{ c_0' } \leq \norm{A}
\]
implies $\LL_{ 1 , \infty }^{ \KK } = 1$.

\end{remark}
\noindent Before we proceed further, we state an interpolation result that will be used repeatedly (Lemmas~\ref{Laa^*^K=1} and~\ref{cota-La1^C}, and Theorem~\ref{main-Lab^R}): it is a direct consequence of the more general interpolation result from~\cite{BePa61}*{p.~302}.

\begin{proposition} \label{interpol}

Let \( a_0 , a_1 , b_0 , b_1 \in [ 1 , \infty ] \), let \( \theta \in ( 0 , 1 ) \), and let \( A \colon c_0 \times c_0 \to \KK \) be a continuous bilinear form such that \( \norm{A}_{ a_0 , b_0 } < \infty \) and \( \norm{A}_{ a_1 , b_1 } < \infty \). If \( a , b \in [ 1 , \infty ] \) are given by
\begin{align*}
\frac{1}{a} & = \frac{ \theta }{ a_0 } + \frac{ 1 - \theta }{ a_1 } \, ; \\[2mm]
\frac{1}{b} & = \frac{ \theta }{ b_0 } + \frac{ 1 - \theta }{ b_1 } \, ,
\end{align*}
then we have
\[
\Aab \leq \bigl ( \norm{A}_{ a_0 , b_0 } \bigr )^{ \theta } \cdot \bigl ( \norm{A}_{ a_1 , b_1 } \bigr )^{ 1 - \theta } \, .
\]

\end{proposition}

\pagebreak

\begin{lemma} \label{Laa^*^K=1}

For all \( a \in [ 2 , \infty ] \) and every continuous bilinear form \( A \colon c_0 \times c_0 \to \KK \) we have \( \norm{A}_{ a , a^{ \ast } } \leq \norm{A} \).

\end{lemma}

\begin{proof}

Observe that the extreme case $a = \infty$ is given by Lemma~\ref{La1^K<=2^1/a}. On the other hand, by Remark~\ref{Lab^K=1;ab>=2}, for every continuous bilinear form $A \colon c_0 \times c_0 \to \KK$ we have $\norm{A}_{ 2 , 2 } \leq \norm{A}$, which settles the case $a = 2$. Finally, given $a \in ( 2 , \infty )$, let $\theta_0 = 1 - \frac{2}{a} \in ( 0 , 1 )$; since
\begin{align*}
\frac{1}{a} & = \frac{ \theta_0 }{ \infty } + \frac{ 1 - \theta_0 }{2} \, ; \\[2mm]
\frac{1}{ a^{ \ast } } & = \frac{ \theta_0 }{1} + \frac{ 1 - \theta_0 }{2} \, ,
\end{align*}
the estimates for $a = 2$ and $a = \infty$ together with Proposition~\ref{interpol} yield
\begin{align*}
\norm{A}_{ a , a^{ \ast } } & \leq \bigl ( \norm{A}_{ \infty , 1 } \bigr )^{ \theta_0 } \cdot \bigl ( \norm{A}_{ 2 , 2 } \bigr )^{ 1 - \theta_0 } \\
& \leq \norm{A}^{ \theta_0 } \cdot \norm{A}^{ 1 - \theta_0 } \\
& = \norm{A} \, . \qedhere
\end{align*}

\end{proof}

\begin{figure}[tbh!]
%
% sets the figure centered with respect to the graphic without the caption, which is precisely the width of Figure 3 (which has no caption)
\settowidth{\figwidth}{\includegraphics{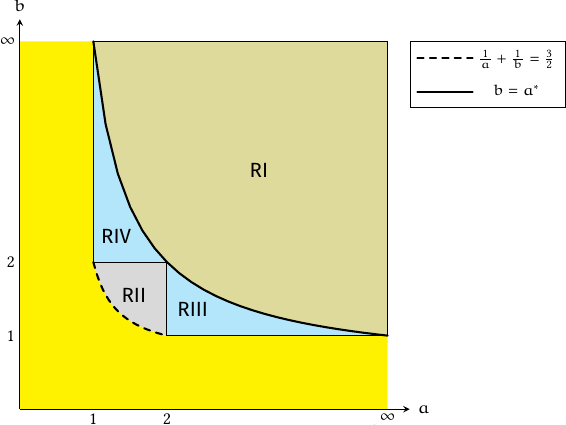}}
\hspace*{\figwidth - \figwidthwithoutcaption}
\includegraphics{Figure2.pdf}
\caption{Illustration of Theorem~\ref{main-Lab^R} and its proof (all the regions are closed): Regions \textsf{RII--RIV} correspond to the value $2^{ \frac{1}{a} + \frac{1}{b} - 1 }$ for the optimal constant, whereas Region~\textsf{RI} corresponds to points for which the optimal constant is $1$. This settles completely the real case.}
\label{fig2}
\end{figure}

\noindent Now we are ready to prove the main result of this work. To this end, we split our reasoning into cases, according to the division of the set of admissible pairs $( a , b )$ into the four subregions determined by the following conditions (see Figure~\ref{fig2}):

\begin{enumerate}[
label   = \textsf{R\Roman*.} ,
align   = left ,
topsep  = 5pt ,
itemsep = 4pt ,
]

\item $a \in [ 1 , \infty ]$ and $b \in [ a^{ \ast } , \infty ]$\,.

\item $a , b \in [ 1 , 2 ]$ and $\frac{1}{a} + \frac{1}{b} \leq \frac{3}{2}$\,.

\item $a \in [2 , \infty ]$ and $b \in [ 1 , a^{ \ast } ]$\,.

\item $a \in [ 1 , 2 ]$ and $b \in [ 2 , a^{ \ast } ]$\,.

\end{enumerate}

\pagebreak

\begin{proof}[Proof of Theorem~\ref{main-Lab^R}]

Let $a , b \in [ 1 , \infty ]$ with $\frac{1}{a} + \frac{1}{b} \leq \frac{3}{2}$\,.

\bigskip

\begin{flushleft}
$\bm{ ( a , b ) \in }$ \textbf{\textsf{RII}}
\end{flushleft}
In this case the result follows from Equation~\eqref{Lab^R-incomplete}.

\bigskip

\begin{flushleft}
$\bm{ ( a , b ) \in }$ \textbf{\textsf{RIII}}
\end{flushleft}
We have $a \in [ 2 , \infty ]$ and $b \in [ 1 , a^{ \ast } ]$. If $b = 1$ or $b = a^{ \ast }$, then $\norm{A}_{ a , b } \leq 2^{ \frac{1}{a} + \frac{1}{b} - 1 } \cdot \norm{A}$ by Lemmas~\ref{La1^K<=2^1/a} and~\ref{Laa^*^K=1}, respectively. If $b \in ( 1 , a^{ \ast } )$, then $b^{ \ast } \in ( a , \infty )$, hence $\theta_1 = 1 - \frac{a}{ b^{ \ast } }$ satisfies $\theta_1 \in ( 0 , 1 )$. Moreover, we have
\begin{alignat*}{4}
\frac{1}{a} \cdot \theta_1 & \ \, + & \frac{1}{a} & \cdot ( 1 - \theta_1 ) && = \frac{1}{a}
\intertext{and}
\frac{1}{1} \cdot \theta_1 & \ \, + \ & \frac{1}{ a^{ \ast } } & \cdot ( 1 - \theta_1 ) && = \frac{1}{1} \cdot \biggl ( 1 - \frac{a}{ b^{ \ast } } \biggr ) + \biggl ( 1 - \frac{1}{a} \biggr ) \cdot \frac{a}{ b^{ \ast } } \\[2mm]
&&&&& = 1 - \frac{1}{ b^{ \ast } } \\[2mm]
&&&&& = \frac{1}{b} \, .
\end{alignat*}
Therefore, by Proposition~\ref{interpol} and Lemmas~\ref{La1^K<=2^1/a} and~\ref{Laa^*^K=1} we get
\begin{align*}
\Aab & \leq \bigl ( \norm{A}_{ a , 1 } \bigr )^{ \theta_1 } \cdot \bigl ( \norm{A}_{ a , a^{ \ast } } \bigr )^{ 1 - \theta_1 } \\
& \leq \bigl ( 2^{ \frac{1}{a} } \cdot \norm{A} \bigr )^{ 1 - \frac{a}{ b^{ \ast } } } \cdot \norm{A}^{ \frac{a}{ b^{ \ast } } } \\
& = \bigl ( 2^{ \frac{1}{a} } \bigr )^{ 1 - \frac{a}{ b^{ \ast } } } \cdot \norm{A} \\
& = 2^{ \frac{1}{a} + \frac{1}{b} - 1 } \cdot \norm{A} \, .
\end{align*}
In addition, the bilinear form $A_0$ given by
\[
A_0 ( \xx \, , \yy ) = x_1 y_1 + x_1 y_2 + x_2 y_1 - x_2 y_2
\]
satisfies $\norm{ A_0 } = 2$ and
\[
\norm{ A_0 }_{ a , b } = \Biggl ( \, \sum_{ k = 1 }^2 \biggl ( \, \sum_{ j = 1 }^2 \, \abs[\big]{ A_0 \bigl ( \ee^{ (k) } , \ee^{ (j) } \bigr ) }^a \biggr )^{ \frac{b}{a} } \Biggr )^{ \frac{1}{b} } = 2^{ \frac{1}{b} + \frac{1}{a} - 1 } \cdot \norm{ A_0 }
\]
for all $a , b \in ( 0 , \infty ]$ (not just for the values $a$ and $b$ considered in this case; see~\cite{diniz}), which shows that the value $2^{ \frac{1}{b} + \frac{1}{a} - 1 }$ cannot be decreased. Thus, in this case we have $\LL_{ a , b }^{ \RR } = 2^{ \frac{1}{b} + \frac{1}{a} - 1 }$.

\bigskip

\begin{flushleft}
$\bm{ ( a , b ) \in }$ \textbf{\textsf{RIV}}
\end{flushleft}
We have $a \in [ 1 , 2 ]$ and $b \in [ 2 , a^{ \ast } ]$. If $b = \infty$, then $a = 1$ (because $b \leq a^{ \ast }$), and the result follows from Remark~\ref{L1infty=Linfty1=1} ($\LL_{ 1 , \infty }^{ \RR } = 1 = 2^{ \frac{1}{1} + \frac{1}{ \infty } - 1 }$). If $b \neq \infty$, then $( b , a )$ lies in Region~\textsf{RIII}, hence $\LL_{ b , a }^{ \RR } = 2^{ \frac{1}{b} + \frac{1}{a} - 1 }$, and since $0 < a \leq b < \infty $, Remark~\ref{minkowski} implies $\LL_{ a , b }^{ \RR } \leq \LL_{ b , a }^{ \RR }$\,. Therefore $\LL_{ a , b }^{ \RR } \leq 2^{ \frac{1}{a} + \frac{1}{b} - 1 }$, and the opposite inequality follows by considering the same bilinear form $A_0$ from the Region~\textsf{RIII} case.

\bigskip

\begin{flushleft}
$\bm{ ( a , b ) \in }$ \textbf{\textsf{RI}}
\end{flushleft}
We have $a \in [ 1 , \infty ]$ and $b \in [ a^{ \ast } , \infty ]$, hence $\Aab \leq \norm{A}_{ a , a^{ \ast } }$ by the embedding of the $\ell_p$ spaces. Also, since the point $( a , a^{ \ast } )$ lies in one of the Regions \textsf{RIII} or \textsf{RIV}, it follows that $\norm{A}_{ a , a^{ \ast } } \leq 2^{ \frac{1}{a} + \frac{1}{ a^{ \ast } } - 1 } \cdot \norm{A} = \norm{A}$. Consequently we have $\Aab \leq \norm{A}$, which shows that $\LL_{ a , b }^{ \RR } = 1$.
\end{proof}

\subsection{The complex version} \label{subsection-complex}

\noindent In this subsection we partially replicate, for the case of complex scalars, the results previously obtained. First, we have the following improvement of Lemma~\ref{La1^K<=2^1/a} in the complex setting:

\begin{lemma} \label{cota-La1^C}

For all \( a \in [ 2 , \infty ] \) and every continuous bilinear form \( A \colon c_0 \times c_0 \to \CC \) we have
\[
\norm{A}_{ a , 1 } \leq \biggl ( \frac{4}{ \pi } \biggr )^{ \frac{1}{a} } \cdot \norm{A} \, .
\]

\end{lemma}

\begin{proof}

By Equation~\eqref{L_12^C;L_21^C} and Proposition~\ref{radem-cota-2^1/a} we have, respectively:
\begin{align*}
\norm{A}_{ 2 , 1 } & \leq \biggl ( \frac{4}{ \pi } \biggr )^{ \frac{1}{2} } \cdot \norm{A} \, ; \\[2mm]
\norm{A}_{ \infty , 1 } & \leq \norm{A} \, ,
\end{align*}
which settles the cases $a = 2$ and $a = \infty$. If $a \in ( 2 , \infty )$, then $\theta = \frac{2}{a} \in ( 0 , 1 )$ satisfies
\begin{align*}
\frac{1}{a} & = \frac{ \theta }{ 2 } + \frac{ 1 - \theta }{ \infty } \, ; \\[2mm]
\quad \frac{1}{1} & = \frac{ \theta }{1} + \frac{ 1 - \theta }{1} \, ,
\end{align*}
and therefore the estimates for $a = 2$ and $a = \infty$ together with Proposition~\ref{interpol} yield
\begin{align*}
\norm{A}_{ a , 1 } & \leq \bigl ( \norm{A}_{ 2 , 1 } \bigr )^{ \theta } \cdot \bigl ( \norm{A}_{ \infty , 1 } \bigr )^{ 1 - \theta } \\
& \leq \biggl ( \frac{4}{ \pi } \biggr )^{ \frac{ \theta}{2} } \cdot \norm{A} \\[2mm]
& = \biggl ( \frac{4}{ \pi } \biggr )^{ \frac{1}{a} } \cdot \norm{A} \, . \qedhere
\end{align*}

\end{proof}

\begin{figure}[tbh!]
\includegraphics{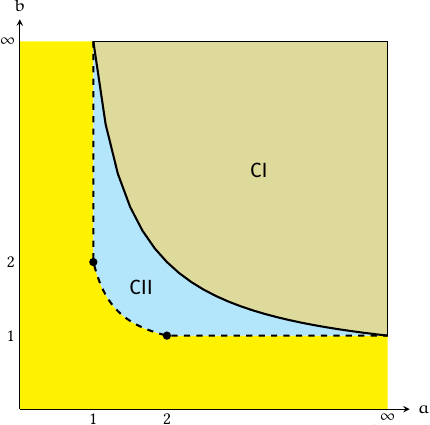}
\caption{Region~\textsf{CI} (closed) consist of points for which the optimal constant in the complex case is $1$ (Theorem~\ref{main-Lab^C}). Region~\textsf{CII} (including the dashed curve minus the points $( 1 , \infty )$ and $( \infty , 1 )$) corresponds to points for which the optimal constant in the complex case is unknown, except for the cases $( a , b ) = ( 1 , 2 )$ and $( a , b ) = ( 2 , 1 )$, for which we have $\LL_{ a , b }^{ \CC } = \frac{2}{ \sqrt{ \pi } }$. (Theorem~\ref{main-Lab^C} provides an upper bound for such unknown values.)}
\label{fig3}
\end{figure}

\noindent The complex partial counterpart of Theorem~\ref{main-Lab^R} can be proved by mimicking its proof, just restricting our attention to the upper estimates, with the following modification: we use Equation~\eqref{Lab^C} instead of Equation~\eqref{Lab^R-incomplete}, and we invoke Lemma~\ref{cota-La1^C} instead of Lemma~\ref{La1^K<=2^1/a}. In this way we obtain the following result (see Figure~\ref{fig3}):

\pagebreak

\begin{theorem} \label{main-Lab^C}

For any \( a , b \in [ 1 , \infty ] \) with \( \frac{1}{a} + \frac{1}{b} \leq \frac{3}{2} \) we have
\[
\LL_{ a , b }^{ \CC } \leq \biggl ( \frac{4}{ \pi } \biggr )^{ \max \{ 0 , \frac{1}{a} + \frac{1}{b} - 1 \} } \, .
\]
In particular, for \( a , b \in [ 1 , \infty ] \) with \( \frac{1}{a} + \frac{1}{b} \leq 1 \) we have \( \LL_{ a , b }^{ \CC } = 1 \).

\end{theorem}

\section{Variants of \texorpdfstring{\noun{Khinchin}}{Khinchin}'s inequality}

\noindent \noun{Khinchin}'s inequality for \noun{Steinhaus} variables plays a crucial role in the improvement of estimates for the constants in inequalities of the type \noun{Littlewood}'s $4 / 3$ (and their extended multilinear versions) for complex scalars (see~\cites{ABPS, cnpr}). For a very recent approach to \noun{Khinchin}'s inequality we refer to~\cite{rama}.

In this section we obtain, among other results, a version of \noun{Khinchin}'s inequality for \noun{Steinhaus} variables (Theorem~\ref{Sr=L1r^C=Lr1^C}). This result is an analog of Proposition~\ref{radem-cota-2^1/a}, which deals with \noun{Rademacher} functions; see Equation~\eqref{rademacher}. The proof of Proposition~\ref{radem-cota-2^1/a} depends crucially on the Contraction Principle from~\cite{diestel}*{p.~231}, which is stated only for randomized sums (independent symmetric real-valued random variables), and we do not know if a similar principle holds for independent symmetric complex-valued random variables.

For this reason, our proof of Theorem~\ref{Sr=L1r^C=Lr1^C} relies instead on ideas from~\cite{cnpr}, where another kind of \noun{Khinchin}'s inequality previously developed by \noun{Ron Blei} is studied (we also consider a variant of such inequality; see Theorem~\ref{blei-khinchin}). In this way, we prove that a certain optimal constant $\St_r$ (see Equation~\ref{def-Sr}) coincides with the optimal constants $\LL_{ 1 , r }^{ \CC }$ and $\LL_{ r, 1 }^{ \CC }$, for all $r \in [ 2 , \infty ]$; this is precisely the analog of Proposition~\ref{radem-cota-2^1/a}: $2^{ \frac{1}{r} } = \LL_{ 1 , r }^{ \RR } = \LL_{ r , 1 }^{ \RR }$ (Theorem~\ref{main-Lab^R}).

\pagebreak

\noindent \emph{From now on, \( N \) denotes a positive integer and \( \aa \) denotes a tuple \( ( a_1 , \dotsc , a_N ) \in \CC^N \)}.

\bigskip

\noindent It is well known that, by considering the dyadic expansion of $2^N t$ for each $t \in [ 0 , 1 ]$, one can prove the equality
\[
\int_0^1 \ \abs[\Bigg]{ \sum_{ j = 1 }^N r_j (t) \, a_j } \, \mathrm{d} t = \frac{1}{ 2^N } \sum_{ \tt \in \{ \pm 1 \}^N } \abs[\Bigg]{ \sum_{ j = 1 }^N a_j t_j } \, ,
\]
where the functions $r_j$ are given by Equation~\eqref{rademacher}. In this way, Proposition~\ref{radem-cota-2^1/a} can be restated as follows: for all $r \in [ 2 , \infty ]$ and any $\aa$ we have
\begin{equation} \label{eq-E_2(a)}
\biggl ( \, \sum_{ j = 1 }^N \, \abs{ a_j }^r \biggr )^{ \frac{1}{r} } \leq 2^{ \frac{1}{r} }\Biggl ( \frac{1}{ 2^N } \sum_{ \tt \in \{ \pm 1 \}^N } \, \abs[\Bigg]{ \sum_{ j = 1 }^N a_j t_j } \, \Biggr ) \, .
\end{equation}
The counterpart for the average
\[
\frac{1}{ 2^N } \sum_{ \tt \in \{ \pm 1 \}^N } \,\abs[\Bigg]{ \sum_{ j = 1 }^N a_j t_j }
\]
in the complex framework is
\begin{equation} \label{int-E_2(a)}
\biggl ( \frac{1}{ 2 \pi } \biggr )^N \int_0^{ 2 \pi } \dotsi \int_0^{ 2 \pi } \abs[\Bigg]{ \sum_{ j = 1 }^N a_j \mathrm{e}^{ \bm { \mathrm{i} } t_j } } \, \mathrm{d} t_1 \dotsb \mathrm{d} t_N \, .
\end{equation}
Notice that this quantity is precisely the expected value
\begin{equation}
\EE \abs[\Bigg]{ \sum_{ j = 1 }^N a_j \varepsilon_j } \, ,
\end{equation}
where $( \varepsilon_1 , \dotsc , \varepsilon_N )$ is a \emph{tuple of \noun{Steinhaus} variables} (that is, a tuple of independent and identically distributed random variables with uniform distribution on the complex unit circle).

\noun{Jerzy Sawa} showed that $\frac{2}{ \sqrt{ \pi } }$ is the least value $\St$ such that
\[
\Biggl ( \, \sum_{ j = 1 }^N \, \abs{ a_j }^2 \Biggr )^{ \frac{1}{2} } \leq \St \cdot \EE \abs[\Bigg]{ \sum_{ j = 1 }^N a_j \varepsilon_j }
\]
for all $N$ and any $\aa$; see~\cite{sawa}. This result is known as \emph{\noun{Khinchin}'s inequality for \noun{Steinhaus} variables}. A generalization, in the spirit of Proposition~\ref{radem-cota-2^1/a}, follows from the embedding of the $\ell_p$ spaces: for all $r \in [ 2 , \infty ]$ there is an optimal constant $\St_r \in [ 1 , \infty )$ such that
\begin{equation} \label{def-Sr}
\Biggl ( \, \sum_{ j = 1 }^N \, \abs{ a_j }^r \Biggr )^{ \frac{1}{r} } \leq \St_r \cdot \EE \abs[\Bigg]{ \sum_{ j = 1 }^N a_j \varepsilon_j }
\end{equation}
for any $N \geq 1$ and any $\aa$. Moreover, we have $\St_r \leq \frac{2}{ \sqrt{ \pi } }$\,.

\pagebreak

Our goal is to prove the following characterization of the values $\St_r$:

\begin{theorem} \label{Sr=L1r^C=Lr1^C}

For all \( r \in [ 2 , \infty ] \) we have \( \St_r = \LL_{ 1 , r }^{ \CC } = \LL_{ r , 1 }^{ \CC } \)\,.

\end{theorem}

\noindent In what follows we state a \noun{Khinchin}-type inequality that extends inequality~\eqref{eq-E_2(a)} (see Remark~\ref{EM(a),M=2}) and recovers inequality~\eqref{def-Sr} as a limiting case (Corollary~\ref{S_r<=L1r}). First, we need to introduce some notation and results.

For each integer $M \geq 2$, let
\[
T_M = \biggl \{ \exp \biggl ( \frac{ 2 \pi j }{M} \, \bm{ \mathrm{i} } \biggr ) \colon j = 0 , \dotsc , M - 1 \biggr \} \, .
\]
In addition, for any bilinear form $A \colon \CC^{ n_1 } \times \CC^{ n_2 } \to \CC$ we define the norm $\norm{A}_M$ by
\[
\norm{A}_M = \sup \, \bigl \{ \abs{ A( \zz \, , \ww ) } \colon \ww \in T_M^{ n_2 } \ ; \ \abs{ z_k } \leq 1 \ \textnormal{for} \ k = 1 , \dotsc , n_1 \bigr \} \, .
\]

\begin{proposition} [\cite{cnpr}*{Theorem~2.3}] \label{cota-||A||_M}

If \( M \geq 3 \), then
\[
\norm{A}_M \leq \norm{A} \leq \sec \biggl ( \frac{ \pi }{M} \biggr ) \cdot \norm{A}_M
\]
for all bilinear forms \( A \colon \CC^{ n_1 } \times \CC^{ n_2 } \to \CC \).

\end{proposition}

\noindent For a tuple $\aa$, we define
\begin{equation} \label{def-E_M(a)}
E_M ( \aa ) = \frac{1}{ M^N } \sum_{ \tt \in T_M^N } \, \abs[\Bigg]{ \sum_{ j = 1 }^N a_j t_j } \, .
\end{equation}
Using the dominated convergence theorem, it is possible to prove that
\begin{equation} \label{E_M-tends-to-E}
\lim_{ M \to \infty } E_M ( \aa ) = \EE \abs[\Bigg]{ \sum_{ j = 1 }^N a_j \varepsilon_j } \, .
\end{equation}
We need the following auxiliary result:

\begin{lemma} [\cite{cnpr}*{Lemma~2.4}] \label{E_M(a)=E_M(aw)}

If \( M \geq 2 \) and \( \ww \in T_M^N \)\,, then for all \( \aa \) we have
\[
E_M ( \aa ) = E_M \bigl ( ( a_j w_j )_{ j = 1 }^N \bigr ) \, .
\]

\end{lemma}

\noindent Now we are ready to state and prove the announced \noun{Khinchin}-type inequality, valid for $r \in [2, \infty ]$. We emphasize that this result is already known for $r = 2$; see~\cite{blei}*{Chapter II \S 6}).

\begin{theorem}[\noun{Blei--Khinchin} inequality] \label{blei-khinchin}

Given \( r \in [ 2 , \infty ] \) and \( M \geq 3 \) we have, for all \( N \) and any \( \aa \):
\[
\Biggl ( \, \sum_{ j = 1 }^N \, \abs{ a_j }^r \Biggr )^{ 1 / r } \leq \LL_{ 1 , r }^{ \CC } \cdot \sec \biggl ( \frac{ \pi }{M} \biggr ) \cdot E_M ( \aa ) \, .
\]

\end{theorem}

\begin{remark} \label{EM(a),M=2}

By Theorem~\ref{main-Lab^R} and Equation~\eqref{def-E_M(a)} we may rewrite inequality~\eqref{eq-E_2(a)} as
\[
\Biggl ( \, \sum_{ j = 1 }^N \, \abs{ a_j }^r \Biggr )^{ \frac{1}{r} } \leq \LL_{ 1 , r }^{ \RR } \cdot E_2 ( \aa ) \, .
\]
In this sense, the result of Theorem~\ref{blei-khinchin} constitutes a complement for inequality~\eqref{eq-E_2(a)}.

\end{remark}

\begin{proof}[Proof of Theorem~\ref{blei-khinchin}]

Let us fix an enumeration of the set $T_M^N$, say $\bm{ \tau }^{ (1) } , \dotsc , \bm{ \tau }^{ ( M^N ) }$, with $\bm{ \tau }^{ (k) } = \bigl ( \tau_1^{ (k) } , \dotsc , \tau_N^{ (k) } \bigr )$ for each $k$. Consider the bilinear form $A \colon \CC^{ M^N } \times \CC^N \to \CC$ given by
\[
A \bigl ( \ee^{ (k) } , \ee^{ (j) } \bigr ) = \frac{ a_j \tau_j^{ (k) } }{ M^N } \quad ( k \in \{ 1 , \dotsc , M^N \} \ ; \ j \in \{ 1 , \dotsc , N \} ) \, .
\]
Let $j \in \{ 1 , \dotsc , N \}$. Since $\abs[\big]{ \tau_j^{ (k) } } = 1$ for each $k$, it follows that
\begin{align*}
\sum_{ k = 1 }^{ M^N } \, \abs[\big]{ A \bigl ( \ee^{ (k) } , \ee^{ (j) } \bigr ) } & = \sum_{ k = 1 }^{ M^N } \frac{ \abs{ a_j } \cdot \abs[\big]{ \tau_j^{ (k) } } }{ M^N } \\[2mm]
& = \abs{ a_j } \, \sum_{ k = 1 }^{ M^N } \frac{1}{ M^N } \\[2mm]
& = \abs{ a_j } \, .
\end{align*}
Thus, by Theorem~\ref{main-Lab^C} and Proposition~\ref{cota-||A||_M} we have
\begin{align*}
\Biggl ( \, \sum_{ j = 1 }^N \,\abs{ a_j }^r \Biggr )^{ \frac{1}{r} } & = \Biggl ( \, \sum_{ j = 1 }^N \biggl ( \, \sum_{ k = 1 }^{ M^N } \, \abs[\big]{ A \bigl ( \ee^{ (k) } , \ee^{ (j) } \bigr ) } \biggr )^r \, \Biggr )^{ \frac{1}{r} } \\[2mm]
& \leq \LL_{ 1 , r }^{ \CC } \cdot \norm{A} \\
& \leq \LL_{ 1 , r }^{ \CC } \cdot \sec \biggl ( \frac{ \pi }{M} \biggr ) \cdot \norm{A}_M \, .
\end{align*}
In order to finish the proof it is sufficient to prove that $\norm{A}_M \leq E_M ( \aa )$, which amounts to have $\abs{ A( \zz \, , \ww ) } \leq E_M ( \aa )$ for all $\zz \in \CC^{ M^N }$ with $\sup_{ 1 \leq k \leq M^N } \, \abs{ z_k } \leq 1$ and all $\ww \in T_M^N$. But this is true: in fact, for any such $\zz$ and $\ww$ we have
\begin{align*}
\abs{ A( \zz , \ww ) } & = \abs[\Bigg]{ \sum_{ k = 1 }^{ M^N } \sum_{ j = 1 }^N A \bigl ( \ee^{ (k) } , \ee^{ (j) } \bigr ) z_k w_j } \\[2mm]
& = \frac{1}{ M^N } \abs[\Bigg]{ \sum_{ k = 1 }^{ M^N } \Biggl ( z_k \cdot \sum_{ j = 1 }^N a_j w_j \tau_j^{ (k) } \Biggr ) } \\[2mm]
& \leq \frac{1}{ M^N } \sum_{ k = 1 }^{ M^N } \, \abs{ z_k } \cdot \abs[\Bigg]{\sum_{ j = 1 }^N a_j w_j \tau_j^{ (k) } } \\[2mm]
& \leq \frac{1}{ M^N } \sum_{ k = 1 }^{ M^N } \ \abs[\Bigg]{\sum_{ j = 1 }^N a_j w_j \tau_j^{ (k) } } \, ,
\end{align*}
and the definition of the elements $\tau_j^{ (k) }$ implies that the last term above is equal to
\begin{align*}
\frac{1}{ M^N } \sum_{ \tt \in T_M^N } \ \abs[\Bigg]{ \sum_{ j = 1 }^N a_j t_j w_j } & = E_M \bigl ( ( a_j w_j )_{ j = 1 }^N \bigr ) \\[2mm]
& = E_M ( \aa ) \, ,
\end{align*}
by Lemma~\ref{E_M(a)=E_M(aw)}.
\end{proof}

\begin{corollary} \label{S_r<=L1r}

Given \( r \in [ 2 , \infty ] \) we have, for all \( N \) and any \( \aa \):
\[
\Biggl ( \, \sum_{ j = 1 }^N \, \abs{ a_j }^r \Biggr )^{ \frac{1}{r} } \leq \LL_{ 1 , r }^{ \CC } \cdot \EE \abs[\Bigg]{ \sum_{ j = 1 }^N a_j \varepsilon_j } \, .
\]

\end{corollary}

\begin{proof}

Just take $M \to \infty$ in Theorem~\ref{blei-khinchin}, and use Equation~\eqref{E_M-tends-to-E}.
\end{proof}

\begin{lemma} \label{L_r1<=S_r}

For all \( r \in [ 2 , \infty ] \) and every continuous bilinear form \( A \colon c_0 \times c_0 \to \CC \) we have
\[
\sum_{ k = 1 }^{ \infty } \Biggl ( \, \sum_{ j = 1 }^{ \infty } \, \abs[\big]{ A \bigl ( \ee^{ (k) } , \ee^{ (j) } \bigr ) }^r \Biggr )^{ \frac{1}{r} } \leq \St_r \cdot \norm{A} \, .
\]

\end{lemma}

\begin{proof}

Proceed exactly as in the proof of Lemma~\ref{La1^K<=2^1/a}, by using~\eqref{def-Sr} instead of~\eqref{eq-radem-cota-2^1/a} and $\St_r$ instead of $2^{ \frac{1}{r} }$, and by substituting the integral $\int_0^1 \abs[\big]{ \sum_{ j = 1 }^N r_j (t) \, A \bigl ( \ee^{ (k) } , \ee^{ (j) } \bigr ) } \, \mathrm{d} t$ with the integral formula for the value $\EE \abs[\big]{ \sum_{ j = 1 }^N A \bigl ( \ee^{ (k) } , \ee^{ (j) } \bigr ) \, \varepsilon_j }$, via Equation~\eqref{int-E_2(a)}.
\end{proof}

\noindent Now we can easily prove the main result of this section.

\begin{proof}[Proof of Theorem~\ref{Sr=L1r^C=Lr1^C}]

Let $r \in [ 2 , \infty ]$. By Corollary~\ref{S_r<=L1r} we have $\St_r \leq \LL_{ 1 , r }^{ \CC }$\,, and Lemma~\ref{L_r1<=S_r} implies $\LL_{ r , 1 }^{ \CC }\leq \St_r$\,. Finally, if $r \neq \infty$, then Remark~\ref{minkowski} implies $\LL_{ 1 , r }^{ \CC } \leq \LL_{ r , 1 }^{ \CC }$ (because $r > 1$), whereas for $r = \infty$ we have $\LL_{ 1 , r }^{ \CC } = \LL_{ r , 1 }^{ \CC } = 1$ by Theorem~\ref{main-Lab^C}. Putting together these inequalities we obtain $\St_r = \LL_{ 1 , r }^{ \CC } = \LL_{ r , 1 }^{ \CC }$\,.
\end{proof}

\section{\texorpdfstring{$( q , s )$-cotype constants of $\ell_1$}{(q,s)-cotype constants of \unichar{"1D4C1}\unichar{"2081}}}

\noindent Let $q \in [ 2 , \infty ]$ and $s \in ( 0 , \infty )$. A (real or complex) \noun{Banach} space $X$ is said to have \( ( q , s ) \)-\emph{cotype} (see~\cite{diestel}*{Remark 11.5(a)}) if there is a constant $\Cotype \in [ 1 , \infty )$ such that, for all $N \geq 1$ and any $x_1 , \dotsc , x_N \in X$,
\begin{equation} \label{eq-cotype-(q,s)}
\Biggl ( \, \sum_{ k = 1 }^N \, \norm{ x_k }^q \Biggr )^{ \frac{1}{q} } \leq \Cotype \cdot \Biggl ( \int_0^1 \norm[\Bigg]{ \sum_{ k = 1 }^N r_k (t) \, \cdot x_k }^s \, \mathrm{d} t \Biggr )^{ \frac{1}{s} } \, ,
\end{equation}
where the functions $r_k$ are given by Equation~\eqref{rademacher}. The smallest value of $\Cotype$ is denoted by $\Cotype_{ q , s } (X)$, the \( ( q , s ) \)-\emph{cotype constant of} \( X \).

\pagebreak

The notion of cotype emerged in the early 1970s; however, previous results in disguise about the cotype of $L^p$-spaces existed as far back as the 1930s (by \noun{Władysław Orlicz}). It is one of the cornerstones of modern \noun{Banach} space theory and reflects the interplay between geometry and probability in \noun{Banach} spaces.

Regarding the general theory, we have the following results:

\begin{proposition}[\cite{pelsansan}*{Lemma~1}] \label{max-xj<=rademacher}

For every \noun{Banach} space \( X \) and any \( x_1 , \dotsc , x_N \in X \), we have
\[
\max_{ 1 \leq k \leq N } \ \norm{ x_k } \leq \int_0^1 \norm[\Bigg]{ \sum_{ j = 1 }^N r_j (t) \, \cdot x_j } \, \mathrm{d} t \, .
\]
Equivalently, for every \noun{Banach} space \( X \) we have \( \Cotype_{ \infty , 1 } (X) = 1 \).
\end{proposition}

\begin{proposition}[\noun{Kahane--Khinchin}'s inequality~\cite{diestel}*{pp.~211}] \label{Kahane}

For any \( s , s' \in ( 0 , \infty ) \), there exists an optimal constant \( \Kah_{ s , s' } \in [ 1 , \infty ) \) such that
\[
\Biggl ( \int_0^1 \norm[\Bigg]{ \sum_{ j = 1 }^N r_j (t) \, \cdot x_j }^{ s' } \, \mathrm{d} t \Biggr )^{ \frac{1}{ s' } } \leq \Kah_{ s , s' } \times \Biggl ( \int_0^1 \norm[\Bigg]{ \sum_{ j = 1 }^N r_j (t) \, \cdot x_j }^s \, \mathrm{d} t \Biggr )^{ \frac{1}{s} } \, ,
\]
regardless of the choice of a \noun{Banach} space \( X \) and vectors \( x_1 , \dotsc , x_N \in X \).
\medskip

\noindent Consequently, a \noun{Banach} space has \( ( q , s_0 ) \)-cotype for some \( s_0 \in ( 0 , \infty ) \) if and only if has \( ( q , s ) \)-cotype for all \( s \in ( 0 , \infty ) \), and in such a case we have \( \Cotype_{ q , s } ( X ) \leq \Cotype_{ q , s_0 } ( X ) \cdot \Kah_{ s , s_0 } \)\,.

\end{proposition}
\smallskip
The embedding of the $\ell_p$ spaces applied to the left side of Equation~\eqref{eq-cotype-(q,s)} implies that the cotype constants are decreasing in the parameter $q$: if a \noun{Banach} space $X$ has $( q , s )$-cotype, then it has $( \tilde{q} , s )$-cotype for any $\tilde{q} \in ( q , \infty ]$, with $\Cotype_{ \tilde{q} , s } (X) \leq \Cotype_{ q , s } (X)$.

Concerning the constants $\Kah_{ s , s' }$\,: from the embedding of the spaces $L^p ( [ 0 , 1 ] )$ we get $\Kah_{ s , s' } = 1$ whenever $s' \leq s$; this fact, together with Propositions~\ref{max-xj<=rademacher} and~\ref{Kahane}, implies that $\Cotype_{ \infty , s } (X) = 1$ for every \noun{Banach} space $X$ and all $s \geq 1$.

\begin{proposition}[\cite{latole}*{Remark 2}] \label{Kahane-specific}

For all \( s , s' \) with \( s \in ( 0 , 1 ] \) and \( s' \in [ s , 2 ] \) we have \( \Kah_{ s , s' } = 2^{ \frac{1}{s} - \frac{1}{ s' } } \). As a consequence, for all \( s \in ( 0 , \infty ) \) we have \( \Kah_{ s , 1 } = 2^{ \max \{ \frac{1}{s} - 1 , 0 \} } \).

\end{proposition}
\smallskip
On the other hand, it is well known that the sequence space $\ell_{1} = \ell_1 ( \KK )$ has $( q , s )$-cotype for all $q \in [ 2 , \infty ]$ and all $s \in ( 0 , \infty )$. To the best of our knowledge, the exact value of the constants $\Cotype_{ q , s } ( \ell_{1} )$ is not available in the literature.

\begin{theorem}

For all \( q \in [ 2 , \infty ] \) and \( s \in (0 , \infty ) \) we have \( 2^{ \frac{1}{q} } \leq \Cotype_{ q , s } ( \ell_1 ) \leq 2^{ \frac{1}{q} + \max \{ \frac{1}{s} - 1 , 0 \} } \). In particular, for \( s \in [ 1 , \infty ) \) we have \( \Cotype_{ q , s } ( \ell_1 ) = 2^{ \frac{1}{q} } \).

\end{theorem}

\begin{proof}

Let $q \in [2 , \infty ]$ and $s \in (0 , \infty )$ . Given $N \geq 1$ and vectors $\xx^{ (1) } , \dotsc , \xx^{ (N) } \in \ell_1$, with $\xx^{ (j) } = \bigl ( x_k^{ (j) } \bigr )_{ k \geq 1 }$ for each $j$, we claim that
\begin{equation} \label{auxiliary1}
\Biggl ( \, \sum_{ k = 1 }^N \, \norm[\big]{ \xx^{ (k) } }_{ \ell_1 }^q \Biggr )^{ \frac{1}{q} } \leq \sum_{ j = 1 }^{ \infty } \Biggl ( \, \sum_{ k = 1 }^N \, \abs[\big]{ x_j^{ (k) } }^q \Biggr )^{ \frac{1}{q} } \, . \pagebreak
\end{equation}
In fact, if $q < \infty$, then Remark~\ref{minkowski} with $( a , b ) = ( 1 , q )$ implies
\begin{align*}
\Biggl ( \, \sum_{ k = 1 }^N \, \norm[\big]{ \xx^{ (k) } }_{ \ell_1 }^q \Biggr )^{ \frac{1}{q} } & = \Biggl ( \, \sum_{ k = 1 }^N \biggl ( \, \sum_{ j = 1 }^{ \infty } \, \abs[\big]{ x_j^{ (k) } } \biggr )^q \Biggr )^{ \frac{1}{q} } \\[2mm]
& \leq \sum_{ k = 1 }^{ \infty } \Biggl ( \, \sum_{ j = 1 }^N \, \abs[\big]{ x_k^{ (j) } }^q \Biggr )^{ \frac{1}{q} } = \sum_{ j = 1 }^{ \infty } \Biggl ( \, \sum_{ k = 1 }^N \, \abs[\big]{ x_j^{ (k) } }^q \Biggr )^{ \frac{1}{q} } \,
\intertext{whereas for $q = \infty$ the inequality follows from}
\max_{ 1 \leq k \leq N } \ \norm[\big]{ \xx^{ (k) } }_{ \ell_1 } & = \max_{ 1 \leq k \leq N } \ \sum_{ j = 1 }^{ \infty } \, \abs[\big]{ x_j^{ (k) } } \\[2mm]
& \leq \sum_{ j = 1 }^{ \infty } \ \max_{ 1 \leq k \leq N } \, \abs[\big]{ x_j^{ (k) } } \, .
\end{align*}
On the other hand, by Proposition~\ref{radem-cota-2^1/a} the right hand term of inequality~\eqref{auxiliary1} satisfies (recall that $q \geq 2$)
\begin{align*} \label{auxiliary2}
\sum_{ j = 1 }^{ \infty } \Biggl ( \, \sum_{ k = 1 }^N \, \abs[\big]{ x_j^{ (k) } }^q \Biggr )^{ \frac{1}{q} } & \leq 2^{ \frac{1}{q} } \sum_{ j = 1 }^{ \infty } \int_0^1 \abs[\Bigg]{ \sum_{ k = 1 }^N r_k (t) \, x_j^{ (k) } } \, \mathrm{d} t \\[2mm]
& = 2^{ \frac{1}{q} } \int_0^1 \sum_{ j = 1 }^{ \infty } \,\abs[\Bigg]{ \sum_{ k = 1 }^N r_k (t) \, x_j^{ (k) } } \, \mathrm{d} t \\[2mm]
& = 2^{ \frac{1}{q} } \int_0^1 \norm[\Bigg]{ \sum_{ k = 1 }^N r_k (t) \, \cdot \xx^{ (k) } }_{ \ell_1 } \! \mathrm{d} t \, .
\end{align*}
This shows that $\ell_1$ has $( q , 1 )$-cotype, with $\Cotype_{ q , 1 } ( \ell_1 ) \leq 2^{ \frac{1}{q} }$. This inequality, together with Propositions~\ref{Kahane} and~\ref{Kahane-specific}, yields
\[\Cotype_{ q , s } ( \ell_1 ) \leq \Cotype_{ q , 1 } ( \ell_1 ) \cdot \Kah_{ s , 1 } = \Cotype_{ q , 1 } ( \ell_1 ) \cdot 2^{ \max \{ \frac{1}{s} - 1 , 0 \} } \leq 2^{ \frac{1}{q} + \max \{ \frac{1}{s} - 1 , 0 \} } \, .
\]
Recall that by Theorem~\ref{main-Lab^R} we have $\LL_{ 1 , q }^{ \RR } = 2^{ \frac{1}{q} }$. Thus, in order to prove the inequality $2^{ \frac{1}{q} } \leq \Cotype_{ q , s } ( \ell_1 )$, we take any continuous bilinear form $A \colon c_0 \times c_0 \to \RR$, and we must show that $\norm{A}_{ 1 , q } \leq \Cotype_{ q , s } ( \ell_1 ) \cdot \norm{A}$; equivalently, for any $N \geq 1$ we must show that
\begin{equation} \label{2^1/q<=C_q,s}
\Biggl ( \, \sum_{ k = 1 }^N \biggl ( \, \sum_{ j = 1 }^N \, \abs[\big]{ A \bigl ( \ee^{ (k) } , \ee^{ (j) } \bigr ) } \biggr )^q \Biggr )^{ \frac{1}{q} } \leq \Cotype_{ q , s } ( \ell_1 ) \cdot \norm{A} \, . \pagebreak
\end{equation}
Let $A_{ N , \ee } \colon c_0 \to \RR^N \subseteq \ell_1$ be the linear map given by
\[
A_{ N , \ee } ( \xx ) = \bigl ( A \bigl ( \xx \, , \ee^{ (j) } \bigr ) \bigr )_{ j = 1 }^N \, .
\]
We have
\begin{align} \label{almost2^1/q<=C_q,s}
\Biggl ( \, \sum_{ k = 1 }^N \biggl ( \, \sum_{ j = 1 }^N \, \abs[\big]{ A \bigl ( \ee^{ (k) } , \ee^{ (j) } \bigr ) } \biggr )^q \Biggr )^{ \frac{1}{q} } & = \Biggl ( \, \sum_{ k = 1 }^N \, \norm[\big]{ A_{ N , \ee } \bigl ( \ee^{ (k) } \bigr ) }_{ \ell_1 }^q \Biggr )^{ \frac{1}{q} } \notag \\[2mm]
& \leq \Cotype_{ q , s } ( \ell_1 ) \cdot \Biggl ( \int_0^1 \norm[\Bigg]{ \sum_{ j = 1 }^N r_j (t) \, \cdot A_{ N , \ee } \bigl ( \ee^{ (j) } \bigr ) }^s _{ \ell_1 } \! \mathrm{d} t \Biggr )^{ \frac{1}{s} } \notag \\[2mm]
& = \Cotype_{ q , s } ( \ell_1 ) \cdot \Biggl ( \int_0^1 \norm[\Bigg]{ A_{ N , \ee } \Biggl ( \sum_{ j = 1 }^N r_j (t) \, \cdot \ee^{ (j) } \Biggr ) }^s _{ \ell_1 } \! \mathrm{d} t \Biggr )^{ \frac{1}{s} } \, . \notag
\intertext{
Since $\norm[\big]{ \sum_{ j = 1 }^N r_j (w) \cdot \ee^{ (j) } }_{ c_0 } \leq 1$ for each $w \in [ 0 , 1 ]$, it follows that the integral in the second factor above satisfies
\begin{align*}
\int_0^1 \norm[\Bigg]{ A_{ N , \ee } \Biggl ( \sum_{ j = 1 }^N r_j (t) \cdot \ee^{ (j) } \Biggr ) }^s _{ \ell_1 } \mathrm{d} t & \leq \sup_{ w \in [ 0 , 1 ] } \ \norm[\Bigg]{ A_{ N , \ee } \Biggl ( \, \sum_{ j = 1 }^N r_j (w) \cdot \ee^{ (j) } \Biggr ) }_{ \ell_1 }^s \\[2mm]
& \leq \sup_{ w \in [ 0 , 1 ] } \ \Biggl [ \, \norm{ A_{ N , \ee } } \cdot \norm[\Bigg]{ \sum_{ j = 1 }^N r_j (w) \cdot \ee^{ (j) } }_{ c_0 } \, \Biggr ]^s \\[2mm]
& \leq \norm{ A_{ N , \ee } }^s \, ,
\end{align*}
and thus we have
}
\Biggl ( \, \sum_{ k = 1 }^N \biggl ( \, \sum_{ j = 1 }^N \, \abs[\big]{ A \bigl ( \ee^{ (k) } , \ee^{ (j) } \bigr ) } \biggr )^q \Biggr )^{ \frac{1}{q} } & \leq \Cotype_{ q , s } ( \ell_1 ) \cdot \bigl ( \norm{ A_{ N , \ee } }^s \bigr )^{ \frac{1}{s} } \notag \\[2mm]
& = \Cotype_{ q , s } ( \ell_1 ) \cdot \norm{ A_{ N , \ee } } \, .
\end{align}
Finally, we claim that $\norm{ A_{ N , \ee } } \leq \norm{A}$, which together with the upper bound~\eqref{almost2^1/q<=C_q,s} yields~\eqref{2^1/q<=C_q,s}. In fact, for any vector $\xx \in c_0$ with $\norm{ \xx }_{ c_0 } \leq 1$ we have (see the proof of Lemma~\ref{La1^K<=2^1/a})
\begin{align*}
\norm{ A_{ N , \ee } ( \xx ) }_{ \ell_1 } & = \sum_{ j = 1 }^N \, \abs[\big]{ A \bigl ( \xx \, , \ee^{ (j) } \bigr ) } \\[2mm]
& \leq \norm{ A( \xx \, , \bm{ \relbar } \, ) }_{ c_0' } \\[2mm]
& \leq \norm{A} \, . \qedhere
\end{align*}

\end{proof}

\end{document}